\documentclass[twocolumn]{article}
\usepackage{mystyle,epsf,amssymb,amsmath}

\newcommand{\ZZ}{\mathbb Z}

\begin{document}

\twocolumn[
\vspace*{30mm}
\centerline{\LARGE Averaged shelling for quasicrystals}\vspace{3ex}
\centerline{\large Michael Baake$^{1}$, Uwe Grimm$^{2}$, 
Dieter Joseph$^{3}$, Przemys{\l}aw Repetowicz$^{2}$}\vspace{2ex}
\begin{footnotesize}
\centerline{
\begin{minipage}{0.85\textwidth}
${}^{1}${\it Institut f\"{u}r Theoretische Physik, 
             Univ.\ T\"{u}bingen,
             Auf der Morgenstelle 14,
             72076 T\"{u}bingen, Germany}\\
${}^{2}${\it Institut f\"{u}r Physik, 
             Technische Universit\"{a}t, 
             09107 Chemnitz, Germany}\\
${}^{3}${\it Max-Planck-Institut f\"{u}r Physik komplexer Systeme, 
             N\"{o}thnitzer Str.~38,
             01187 Dresden, Germany}\\
\end{minipage}
}\end{footnotesize}
\centerline{\footnotesize \today}\vspace{4ex}
\begin{small}
\hrule\vspace{2ex}
\begin{minipage}{\textwidth}
{\bf Abstract}\vspace{2ex}\\
\hp 
The shelling of crystals is concerned with counting the number of
atoms on spherical shells of a given radius and a fixed centre.  Its
straight-forward generalization to quasicrystals, the so-called
central shelling, leads to non-universal answers. As one way to cope
with this situation, we consider shelling averages over all
quasicrystal points.  We express the averaged shelling numbers in
terms of the autocorrelation coefficients and give explicit results
for the usual suspects, both perfect and random.\vspace{2ex}
\end{minipage}\vspace{2ex}
\hrule
\end{small}\vspace{6ex}
]

\section{Introduction}

\begin{figure}[tb]
\centerline{\epsfxsize=0.6\columnwidth\epsfbox{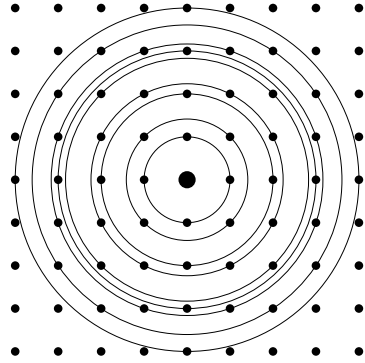}\vspace{-1ex}}
\caption{Shelling of the square lattice. The circles centred at 
the lattice point denoted by the thick dot have radii 
$r^2=1,2,4,5,8,9,10,13,16$, the corresponding shelling numbers
are $\sigma(r)=4,4,4,8,4,4,8,8,4$, compare Eq.~(\protect\ref{eq1}).
\label{fig:square}}
\end{figure}

\hp 
One characteristic feature of a crystal is the number of atoms on
shells of radius $r$ around an arbitrary, but fixed centre, e.g.\
around one fixed atom.  In the simplest case, one thus considers a
lattice, such as the square lattice $\mathbb{Z}^2$ in the plane, and
determines the number of lattice points on circles of radius $r$.  In
the square lattice case, only squared radii $r^2 = m^2 + n^2$ with
$m,n\in\mathbb{Z}$ are possible for obvious geometric reasons, compare
Fig.~\ref{fig:square}.

{}For this example, the answer is well known from number theory \cite{HW}
because it essentially boils down to counting the number of ways that
$r^2$ can be written as the sum of two squares.  The result
reads as follows.  Let $r^2 = M$ be an integer, and $\sigma(r)$ the number
of points of $\mathbb{Z}^2$ on a circle of radius $r$ around the
origin. Then, $\sigma(r) = 4 a(M)$ where $a(M)$ is a multiplicative
arithmetic function, i.e.\ $a(1)=1$ and $a(M N) = a(M) a(N)$ for
coprime $M,N$. So, it suffices to know $a(M)$ for $M$ a prime
power, i.e.\ $M=p^{\ell}$. The explicit result is
\begin{equation}  \label{eq1}
   a(p^{\ell}) \; = \; 
     \begin{cases}
       1      & \text{if $p=2$} \\
       \ell+1 & \text{if $p\equiv 1$ mod 4} \\
       0      & \text{if $p\equiv 3$ mod 4 and $\ell$ odd} \\
       1      & \text{if $p\equiv 3$ mod 4 and $\ell$ even}
     \end{cases}
\end{equation}
There are many explicit formulas known for lattices, see \cite{CS} for
further examples and a systematic exposition by means of lattice theta 
functions.

Now, if one tries to extend this approach to quasicrystals, one
immediate problem arises: due to the lack of periodicity, there are no
``natural'' centres any more, except, perhaps, the centre in a pattern
with exact maximal symmetry --- if this exists! There are several
attempts to find analogues of Eq.~(\ref{eq1}), and explicit results
can be found in \cite{SM,MW,MS,Al,BM}. However, an alternative approach
is desirable in which one does not ask for the shelling numbers around a
fixed centre, but for the {\em averaged} shelling numbers around any
possible centre \cite{BM}, e.g.\ around any centre that is itself in 
the point set. This is what we want to continue to analyze here, in
analogy to earlier work on the coordination numbers of crystals versus
quasicrystals \cite{BG,BGJR}. Note that the ordinary and the averaged
shelling problems coincide for lattices.

\section{A formula for averaged shelling}

\hp 
It is intuitively clear how to define the averaged shelling number for
a point set $\Lambda$: if we fix a radius $r$ and a centred ball
$B_{s}(0)$ of radius $s\gg r$, we inspect any point inside the large
ball, determine the shelling number around it, and average over these
possibilities. This is a radially averaged autocorrelation coefficient.
The problem with it is that there might not be a
limit as $s\to\infty$, and this is indeed a serious problem in
general. However, cut-and-project sets, or model sets as we want to
call them, have the nice property that these averages are guaranteed
to exist. What is more, each finite patch or configuration within a
regular, generic model set (see \cite{Martin} for the precise
conditions) occurs with a uniform frequency. Thus, we can evaluate the
averaged shelling number by means of the existing autocorrelation coefficients.
Let us describe how this works.

Let $\Lambda = \Lambda(\Omega)$ be a regular model set \cite{M} with window 
$\Omega$. The possible distances between points are all of the form $r=|y|$ 
with $y\in\Lambda\!-\!\Lambda$ where
$\Lambda\!-\!\Lambda=\{ x^{}_1 - x^{}_2 \mid x^{}_1, x^{}_2 \in\Lambda\}$. 
We find the formula
\begin{equation}
    \sigma(r) \; = \; \sum_{\stackrel{\scriptstyle y\in\Lambda-\Lambda}
                                     {\scriptstyle |y|=r}} \nu(y) \, .
\end{equation}
Here, for $y\in\Lambda-\Lambda$, $\nu(y)$ is the autocorrelation coefficient,
normalized per point rather than per volume:
\begin{eqnarray}  \label{formula}
    \nu(y) & = &  \lim_{s\to\infty}\frac{1}{|\Lambda_s|}
           \sum_{\stackrel{\scriptstyle x\in\Lambda_s}
                {\scriptstyle x+y\in\Lambda}} 1    \nonumber \\
           & = & \lim_{s\to\infty}\frac{1}{|(\Lambda_s)^*|}
           \sum_{\stackrel{\scriptstyle x^*\in(\Lambda_s)^*}
                {\scriptstyle (x+y)^*\in\Omega}} 1   \\
           & = & \frac{1}{{\rm vol}(\Omega)}
           \int \chi^{}_{\Omega}(x^*)\, \chi^{}_{\Omega}(x^* + y^*)\, dx^*
           \nonumber
\end{eqnarray}
where $\Lambda_s = \Lambda\cap B_s(0)$, $|\Lambda_s|$ is the number of
points in $\Lambda_s$, $(.)^*$ is the $*$-map 
(or lift) from physical to internal space as described in \cite{M},
and  $\chi^{}_{\Omega}$ denotes the characteristic function of the window.
Note that the last step of (\ref{formula}) is correct because Weyl's Theorem 
applies to regular model sets \cite{Martin}.

\section{Result in one dimension}

\hp
Let us explain the typical situation in one dimension with a representative 
example. Consider the silver mean substitution rule $a\mapsto aba$, 
$b\mapsto a$. This gives a semi-infinite fixed point of the form
\begin{equation}
     \mid \!\! abaaabaabaabaaabaabaaabaabaaab \ldots
\end{equation}
If $a$ codes an interval of length $1+\sqrt{2}$ and $b$ one of length 1,
this gives a point set, starting at 0, which is the positive part of
the silver mean model set
\begin{equation}
    \Lambda \; = \; \left\{ \!\!\vphantom{\sqrt{2}}\right.
                    x\in\ZZ[\sqrt{2}\,] \; \mbox{\huge \raisebox{-0.25ex}{$\mid$}}
                    \; x' \in 
                    \Bigl[-\frac{\sqrt{2}}{2},\frac{\sqrt{2}}{2}\,\Bigr]
                    \left.\vphantom{\sqrt{2}}\!\!\right\}
\end{equation}
where $(.)'$ is the algebraic conjugation, defined by 
$\sqrt{2}\mapsto -\sqrt{2}$, and represents the $*$-map in this case.

\begin{figure}[t]
\centerline{\epsfxsize=0.49\columnwidth\epsfbox{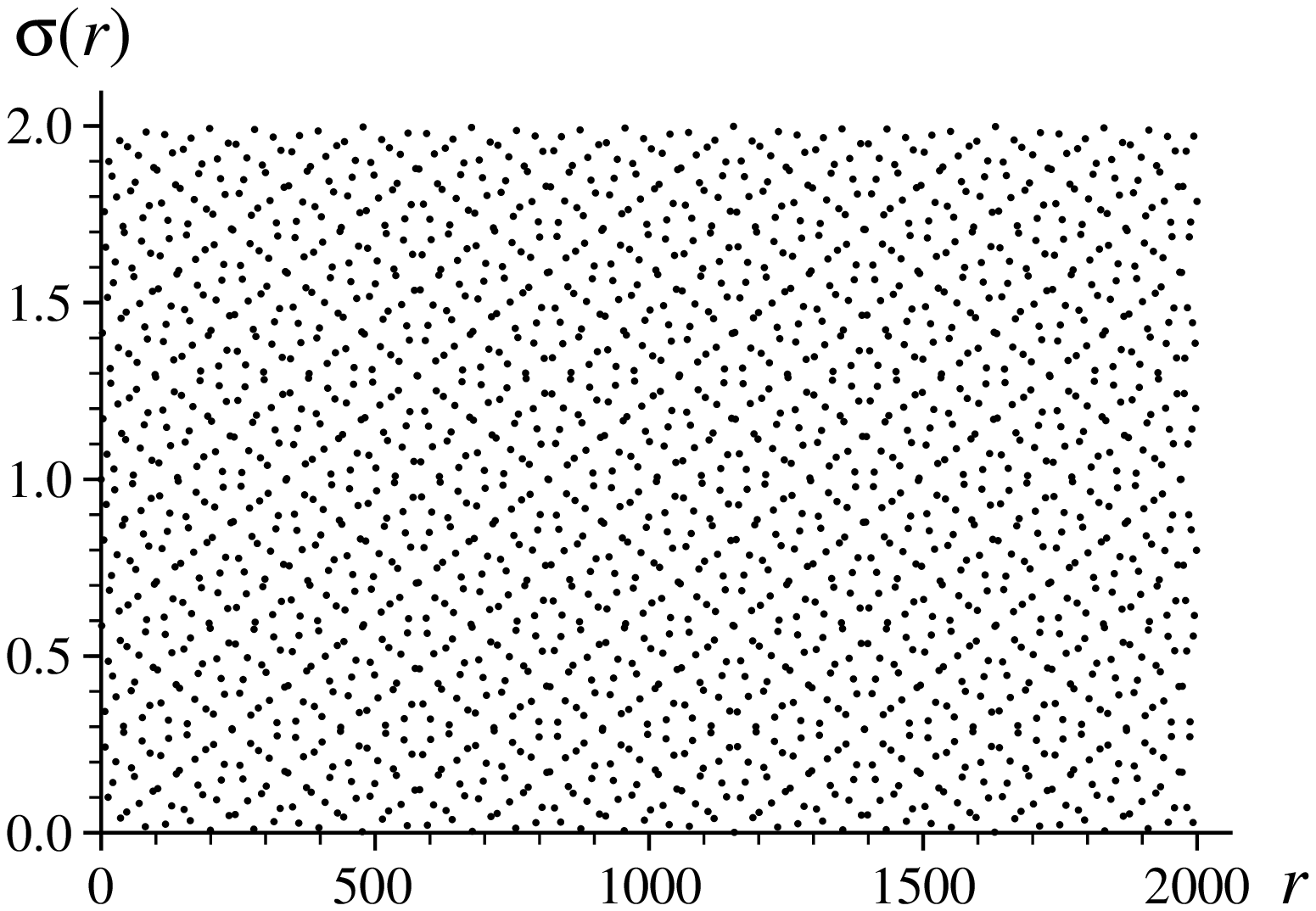}\hfill
\epsfxsize=0.49\columnwidth\epsfbox{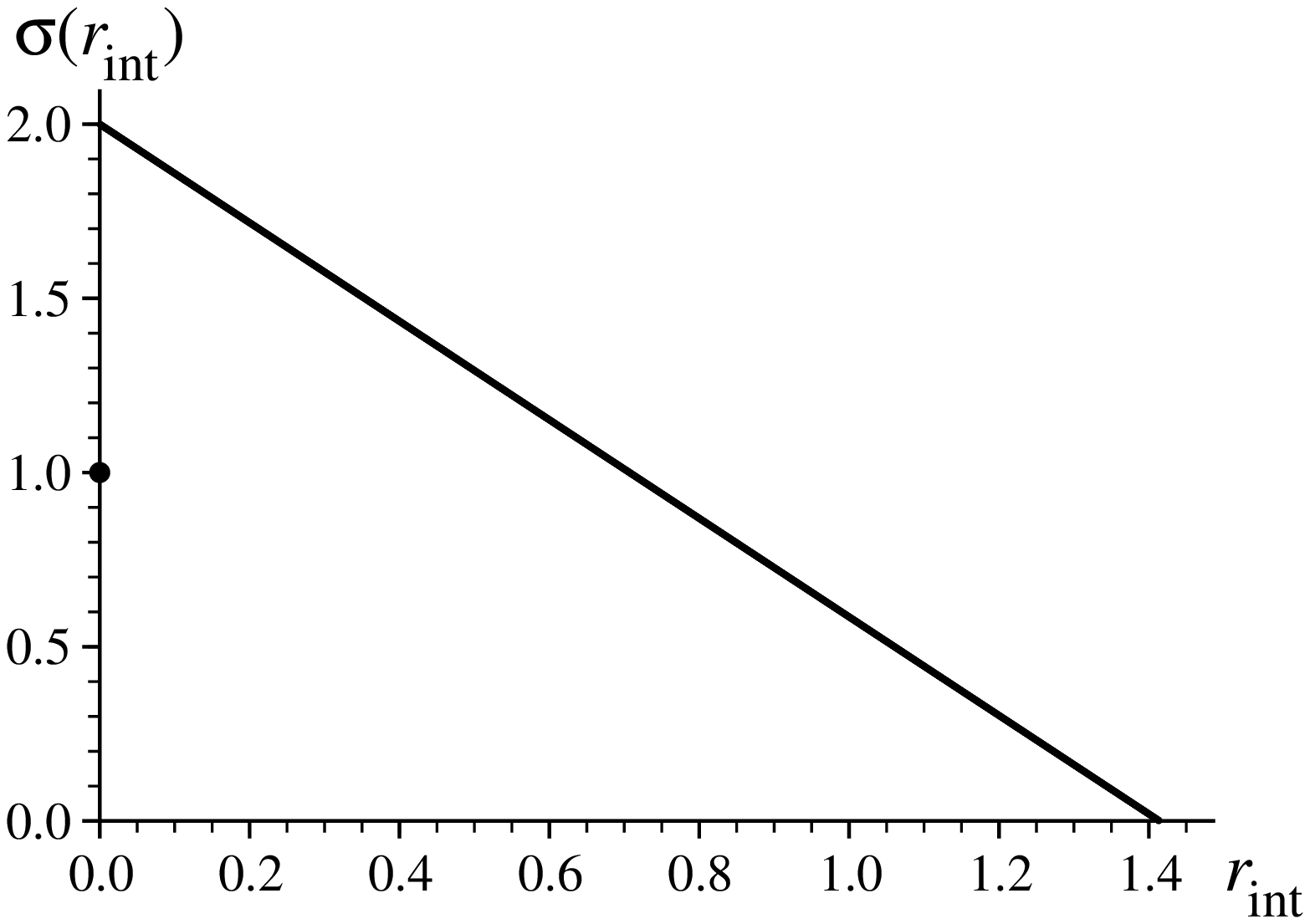}\vspace{-1ex}}
\caption{Averaged shelling function $\sigma(r)$ for the
silver mean model set, plotted against the distance 
$r\in\ZZ[\protect\sqrt{2}]$ in physical space (left) and against 
its counterpart $r^{}_{\rm int}=|r'|$ in internal space 
(right).\label{fig:sm}}
\end{figure}

$\Lambda \!-\! \Lambda$ is a subset of $\ZZ[\sqrt{2}\,]$, and itself
another model set, this time with window $[-\sqrt{2},\sqrt{2}\,]$.
So, applying Eq.~(\ref{formula}), the autocorrelation
coefficient for $y\in\Lambda \!-\! \Lambda$ can be written as $\nu(y) =
f(y')$ with
\begin{equation}  
   f(y') \; = \; 
     \begin{cases}
       1 - \frac{|y'|}{\sqrt{2}}     
              & \text{if $|y'|\leq\sqrt{2}$} \\
       0      & \text{otherwise.} 
     \end{cases}
\end{equation}
The possible shelling radii are then the non-negative numbers in
$\Lambda - \Lambda$.  So, the shelling function is simply
$\sigma(0)=1$ and, since $\nu(y)=\nu(-y)$, $\sigma(y) = 2 f(y')$ for
any positive $y\in\Lambda \!-\! \Lambda$. Now, the averaged shelling can
be calculated, and the smallest positive radii are $1$, $1+\sqrt{2}$,
$2+\sqrt{2}$, $2+2\sqrt{2}$, $3+2\sqrt{2}$ and so on. The result is
plotted in Fig.~\ref{fig:sm}. It is clearly seen 
how the inherent structure of the averaged shelling is hidden in
physical space but clearly visible in internal space.

Let us briefly mention that one can also calculate the averaged shelling
for the corresponding random tiling, as the autocorrelation is explicitly
known in one dimension, see \cite{BM} for details.

\section{Rhombic Penrose tiling}

\begin{figure}[t]
\centerline{\epsfxsize=0.475\columnwidth\epsfbox{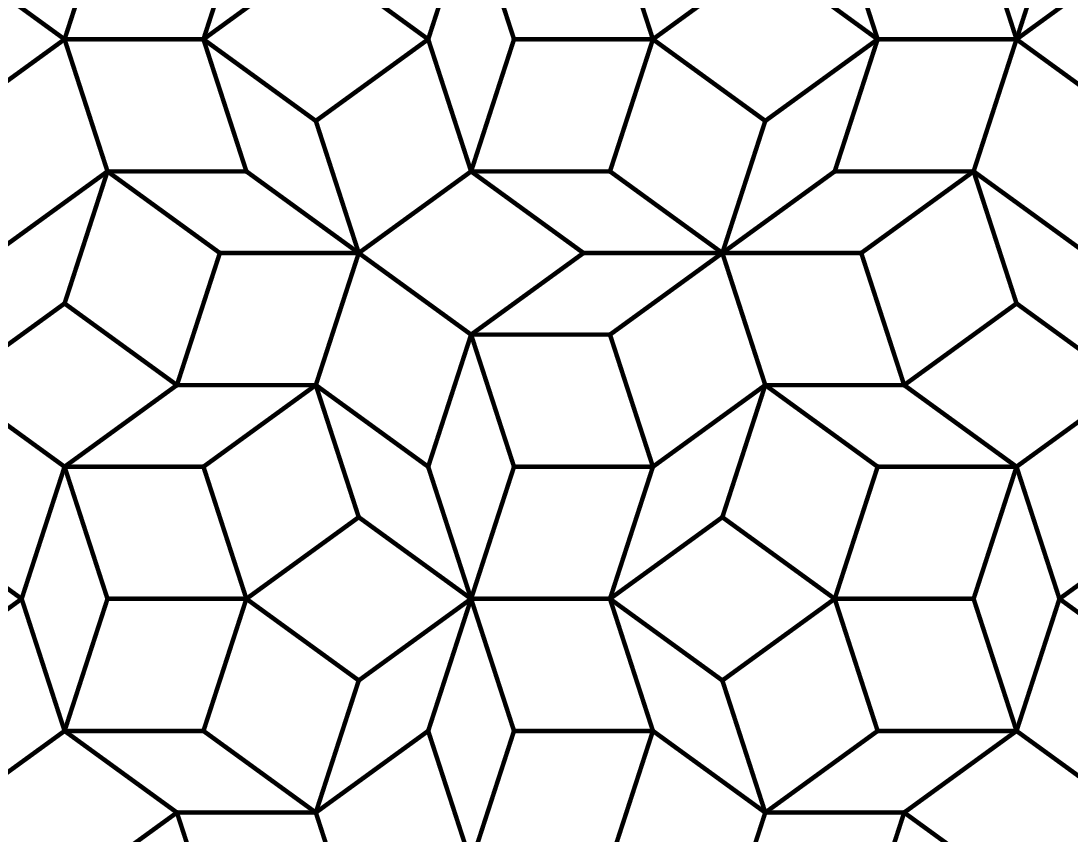}\hfill
\epsfxsize=0.475\columnwidth\epsfbox{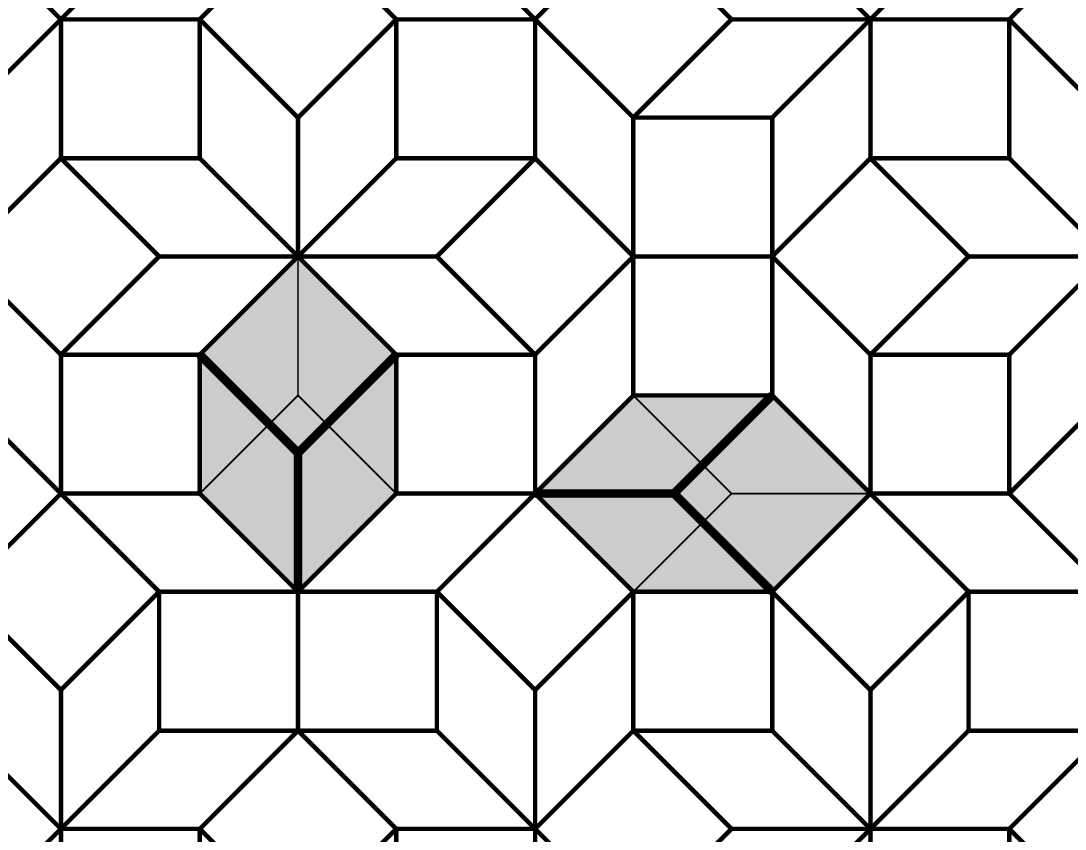}\vspace{-1ex}}
\caption{Patches of the rhombic Penrose (left) and the
Ammann-Beenker tiling (right). For the latter, two simpleton flips
in the shaded hexagons are indicated. The thin lines correspond to
the patch of the perfect tiling. \label{fig:patch}}
\end{figure}

\hp
In the plane, the most studied non-periodic tiling is the classic
Penrose tiling in its rhombic version, see the left part of
Fig.~\ref{fig:patch}.  We consider a realization where the edges have
unit length.  The set of vertex points constitutes a model set,
however one with four different translation classes and hence also
four different windows \cite{BKSZ}.  Although it is still possible to
use Eq.~(\ref{formula}), it is computationally a little easier to
determine all finite patches of a certain radius together with their
frequencies and to calculate the exact averaged shelling from that.  
This is done by a refinement of the window method as described in 
\cite{BKSZ}. The result is summarized in Table~\ref{tab:Pen}.

Note that the squares of possible radii, $r^2$, are totally positive
numbers (i.e.\ both $r^2$ and its algebraic conjugate are positive) in
the ring $\ZZ[\tau]$ where $\tau=(1+\sqrt{5}\,)/2$ is the golden
number. This can be understood from the underlying shelling problem of
the cyclotomic ring $\ZZ[\xi]$ with $\xi = \exp(2\pi i/5)$ and its
number-theoretic ramifications, see \cite{BM} for details and a closed
formula for the central shelling of $\ZZ[\xi]$.

\begin{table}[t]
\caption{Exact results for the first $50$ averaged shelling numbers 
$\sigma(r)$ of the rhombic Penrose
tiling.\protect\rule[-2ex]{0pt}{2ex}\label{tab:Pen}}
\begin{small}
\begin{center}
\begin{tabular}{|r@{\,}c@{\,}r|r@{\,}c@{\,}r||r@{\,}c@{\,}r|r@{\,}c@{\,}r|}
\hline
\multicolumn{3}{|c|}{\rule[-6pt]{0pt}{18pt}$r^2$} &
\multicolumn{3}{c||}{$\sigma(r)$} &
\multicolumn{3}{c|}{$r^2$} &
\multicolumn{3}{c|}{$\sigma(r)$} \\
\hline \rule{0pt}{14pt}
 $2$&$-$&$\tau$   & $4$&$-$&$2\tau$ &
 $5$&$+$&$8\tau$  & $-22$&$+$&$16\tau$ \\
 $1$&   &         & $4$  & & &
$10$&$+$&$5\tau$  & $52$&$-$&$32\tau$ \\
 $3$&$-$&$\tau$   & $8$&$-$&$4\tau$ &
$11$&$+$&$5\tau$  & $-84$&$+$&$52\tau$ \\
 $4$&$-$&$\tau$   & $-12$&$+$&$8\tau$ &
 $8$&$+$&$7\tau$  & $20$&$-$&$12\tau$ \\
 $1$&$+$&$\tau$   & $10$&$-$&$4\tau$ &
 $7$&$+$&$8\tau$  & $-144$&$+$&$92\tau$ \\
 $2$&$+$&$\tau$   & $30$&$-$&$16\tau$ &
$10$&$+$&$7\tau$  & $-106$&$+$&$66\tau$ \\
 $4$&   &         & $10$&$-$&$6\tau$ &
 $7$&$+$&$9\tau$  & $88$&$-$&$52\tau$ \\
 $3$&$+$&$\tau$   & $-28$&$+$&$20\tau$ &
 $8$&$+$&$9\tau$  & $128$&$-$&$76\tau$ \\
 $5$&   &         & $4$&$-$&$2\tau$ &
 $9$&$+$&$9\tau$  & $-138$&$+$&$86\tau$ \\
 $3$&$+$&$2\tau$  & $16$&$-$&$8\tau$ &
 $7$&$+$&$11\tau$ & $150$&$-$&$88\tau$ \\
 $2$&$+$&$3\tau$  & $-4$&$+$&$6\tau$ &
 $9$&$+$&$10\tau$ & $92$&$-$&$56\tau$ \\
 $5$&$+$&$2\tau$  & $-56$&$+$&$36\tau$ &
 $8$&$+$&$11\tau$ & $-12$&$+$&$12\tau$ \\
 $7$&$+$&$\tau$   & $20$&$-$&$12\tau$ &
$11$&$+$&$10\tau$ & $-180$&$+$&$112\tau$ \\
 $3$&$+$&$4\tau$  & $58$&$-$&$32\tau$ &
 $8$&$+$&$12\tau$ & $126$&$-$&$76\tau$ \\
 $5$&$+$&$3\tau$  & $40$&$-$&$24\tau$ &
$13$&$+$&$9\tau$  & $52$&$-$&$32\tau$ \\
 $4$&$+$&$4\tau$  & & & $2\tau$ &
$11$&$+$&$11\tau$ & $156$&$-$&$96\tau$ \\
 $7$&$+$&$3\tau$  & $-64$&$+$&$40\tau$ &
$12$&$+$&$11\tau$ & $92$&$-$&$56\tau$ \\
 $4$&$+$&$5\tau$  & $44$&$-$&$24\tau$ &
$13$&$+$&$11\tau$ & $-32$&$+$&$20\tau$ \\
 $5$&$+$&$5\tau$  & $42$&$-$&$24\tau$ &
$10$&$+$&$13\tau$ & $48$&$-$&$28\tau$ \\
 $7$&$+$&$4\tau$  & $52$&$-$&$32\tau$ &
 $9$&$+$&$14\tau$ & $-80$&$+$&$56\tau$ \\
 $6$&$+$&$5\tau$  & $-68$&$+$&$44\tau$ &
$12$&$+$&$13\tau$ & $-380$&$+$&$236\tau$ \\
 $8$&$+$&$4\tau$  & $20$&$-$&$12\tau$ &
$10$&$+$&$15\tau$ & $122$&$-$&$72\tau$ \\
 $9$&$+$&$4\tau$  & $-32$&$+$&$20\tau$ &
$11$&$+$&$15\tau$ & $-156$&$+$&$100\tau$ \\
 $5$&$+$&$7\tau$  & $-36$&$+$&$28\tau$ &
$13$&$+$&$14\tau$ & $164$&$-$&$100\tau$ \\
 $7$&$+$&$6\tau$  & $120$&$-$&$72\tau$ &
$11$&$+$&$16\tau$ & $108$&$-$&$64\tau$ \\
\hline
\end{tabular}
\end{center}
\end{small}
\vspace*{-3mm}
\end{table}

\section{Octagonal tiling: perfect versus random}

\hp 
Finally, let us consider the vertex set of the well-known
Ammann-Beenker or octagonal tiling as shown in the right half of
Fig.~\ref{fig:patch}.  As above, edges are fixed to have unit
length. The vertex points form again a simple model set, obtained 
(up to scale) from
the four-dimensional lattice $\ZZ^4$ with a regular octagon as window
in the 2D internal space \cite{BJ}.  To obtain a good numerical approximation 
for a larger number of radii in this case, we
have used a large periodic approximant (with $47321$ vertices in the
unit cell) and determined the averages explicitly. For small radii,
the agreement with the exact formula based upon the window technique
is 6 to 8 digits.

\begin{figure}[b]
\centerline{\epsfxsize=0.49\columnwidth\epsfbox{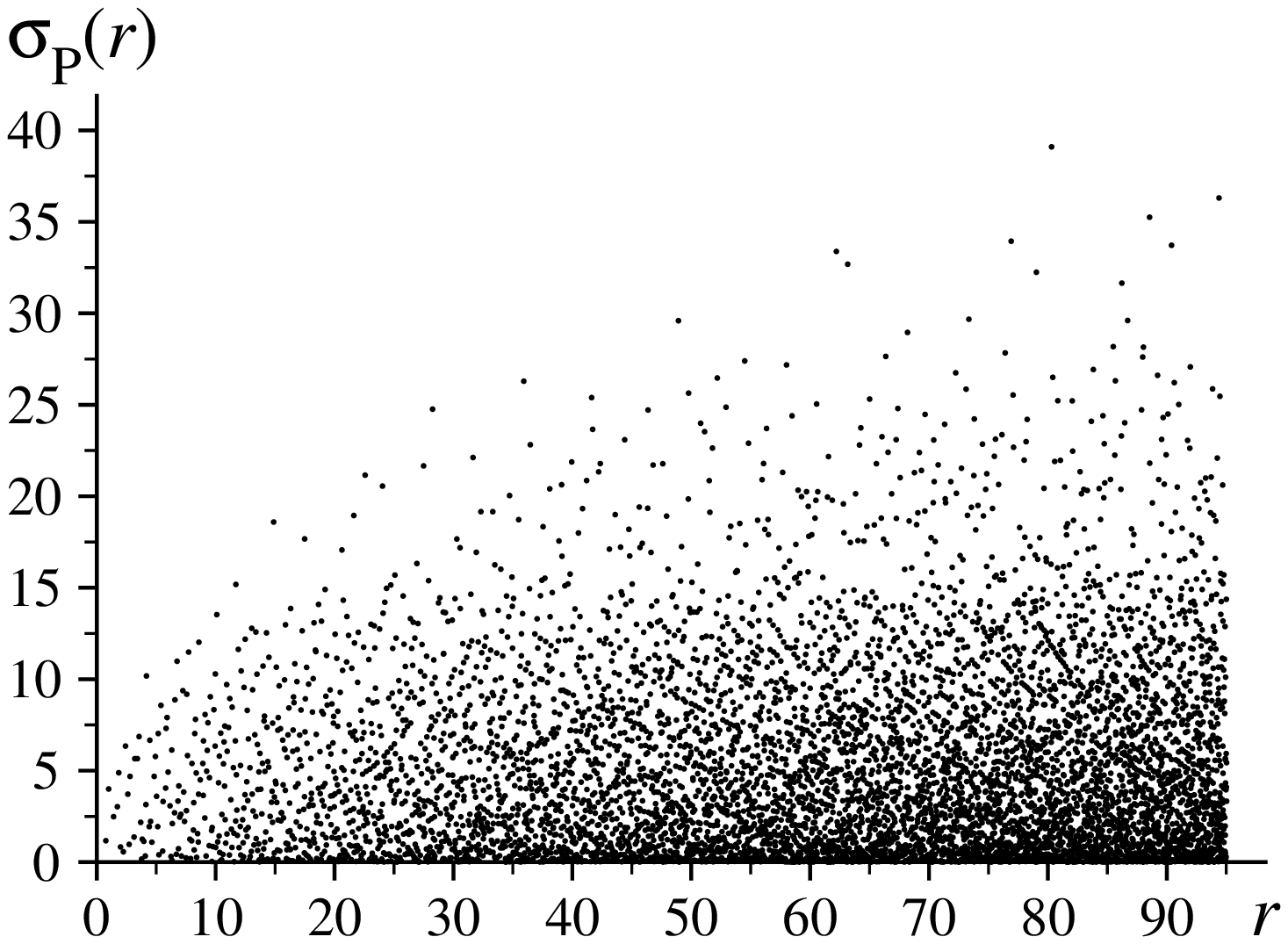}\hfill
\epsfxsize=0.49\columnwidth\epsfbox{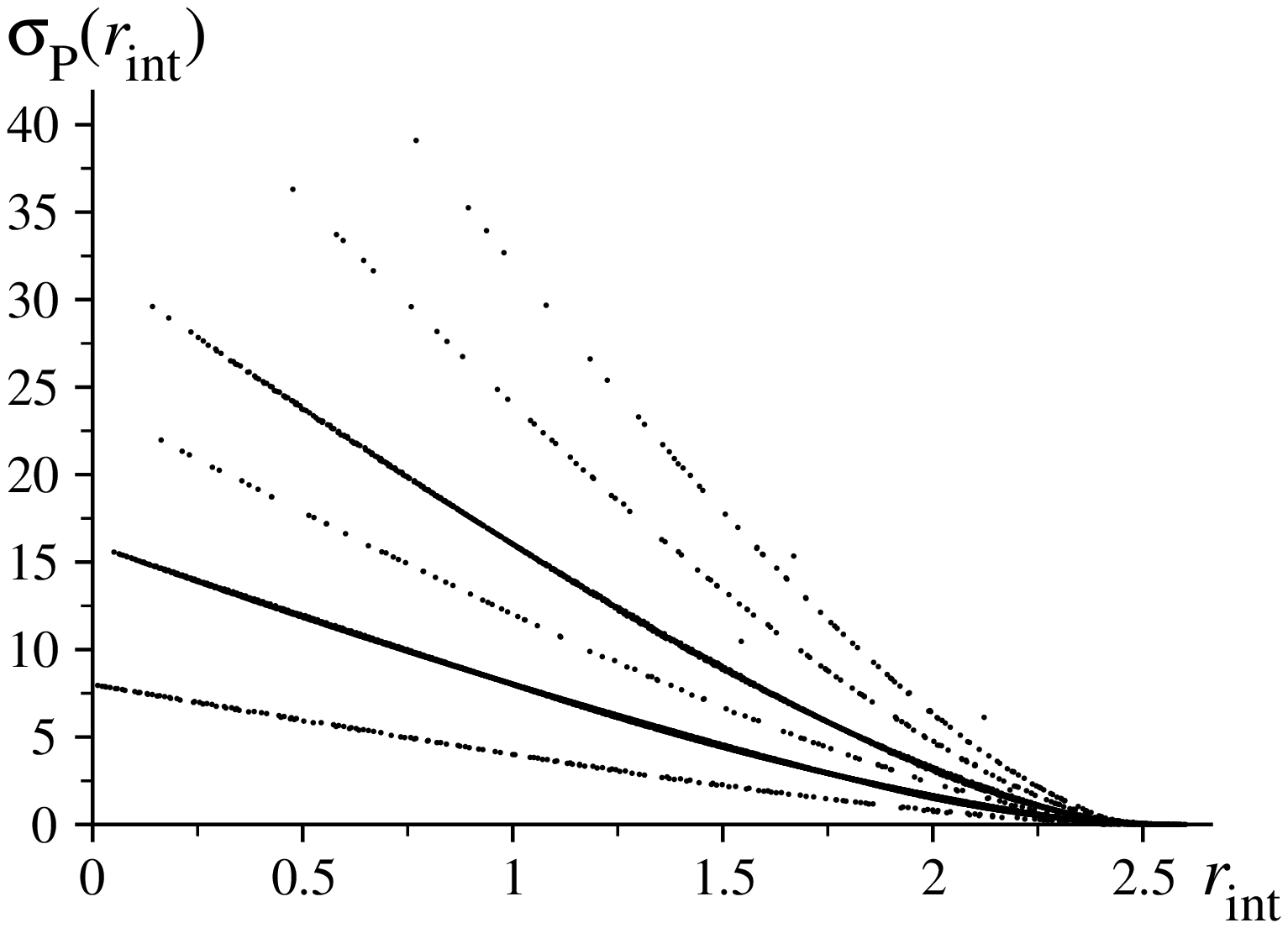}\vspace{-1ex}}
\caption{Averaged shelling function $\sigma_{\rm p}^{}(r)$ for the
perfect Ammann-Beenker tiling, plotted against the radius $r$ in physical 
space (left) and against its counterpart $r^{}_{\rm int}$ in internal space 
(right).\label{fig:AB}}
\end{figure}

\begin{table}[t]
\caption{Averaged shelling function $\sigma_{\rm p}^{}(r)$ for perfect and
$\sigma_{\rm s}^{}(r)$ for stochastic Ammann-Beenker tilings. The possible
shell radii of the perfect tiling are a proper subset of those of the
random tiling.\protect\rule[-2ex]{0pt}{2ex}\label{tab:AB}}
\begin{small}
\begin{center}
\begin{tabular}{|r@{\,}c@{\,}r|r@{\,}c@{\,}r@{$\;$}l|l|}
\hline
\multicolumn{3}{|c|}{\rule[-8pt]{0pt}{20pt}$r^2$} &
\multicolumn{4}{c|}{$\sigma_{\rm p}^{}(r)$} &
\multicolumn{1}{c|}{$\sigma_{\rm s}^{}(r)$} \\
\hline \rule{0pt}{14pt}
 $2$&$-$&$\sqrt{2}$  &   $4$&$-$&$2\sqrt{2} $&$\simeq 1.17157$ & $1.17157$ \\
 $1$&   &            &   $4$&   &            &                 & $4$ \\
 $2$&   &            &  $-6$&$+$&$6\sqrt{2} $&$\simeq 2.48528$ & $2.34807$ \\
 $5$&$-$&$2\sqrt{2}$ &      &   &            &                 & $0.40719$ \\
 $8$&$-$&$4\sqrt{2}$ &      &   &            &                 & $0.01700$ \\
 $3$&   &            &  $20$&$-$&$12\sqrt{2}$&$\simeq 3.02944$ & $2.81463$ \\
 $2$&$+$&$\sqrt{2}$  &  $36$&$-$&$22\sqrt{2}$&$\simeq 4.88730$ & $3.92068$ \\
 $4$&   &            &  $-2$&$+$&$2\sqrt{2} $&$\simeq 0.82843$ & $1.30213$ \\
$13$&$-$&$6\sqrt{2}$ &      &   &            &                 & $0.00246$ \\
 $6$&$-$&$\sqrt{2}$  &      &   &            &                 & $0.57889$ \\
 $5$&   &            & $-56$&$+$&$40\sqrt{2}$&$\simeq 0.56854$ & $1.63616$ \\
$10$&$-$&$3\sqrt{2}$ &      &   &            &                 & $0.04459$ \\
 $3$&$+$&$2\sqrt{2}$ &  $12$&$-$&$4\sqrt{2} $&$\simeq 6.34315$ & $4.11676$ \\
 $6$&   &            &      &   &            &                 & $0.90676$ \\
 $9$&$-$&$2\sqrt{2}$ &      &   &            &                 & $0.09977$ \\
 $4$&$+$&$2\sqrt{2}$ &  $32$&$-$&$20\sqrt{2}$&$\simeq 3.71573$ & $2.73580$ \\
 $6$&$+$&$\sqrt{2}$  & $-16$&$+$&$12\sqrt{2}$&$\simeq 0.97056$ & $1.66129$ \\
 $5$&$+$&$2\sqrt{2}$ &  $16$&$-$&$8\sqrt{2} $&$\simeq 4.68629$ & $4.33600$ \\
 $8$&   &            &      &   &            &                 & $0.29191$ \\
$11$&$-$&$2\sqrt{2}$ &      &   &            &                 & $0.07094$ \\
$14$&$-$&$4\sqrt{2}$ &      &   &            &                 & $0.00739$ \\
 $9$&   &            &      &   &            &                 & $0.62840$ \\
$14$&$-$&$3\sqrt{2}$ &      &   &            &                 & $0.00838$ \\
 $7$&$+$&$2\sqrt{2}$ &  $24$&$-$&$16\sqrt{2}$&$\simeq 1.37258$ & $2.00517$ \\
$10$&   &            &      &   &            &                 & $0.33477$ \\
 $6$&$+$&$3\sqrt{2}$ &      &   &$4\sqrt{2} $&$\simeq 5.65685$ & $4.77325$ \\
$11$ &&              &      &   &            &                 & $0.17416$ \\
$17$&$-$&$4\sqrt{2}$ &      &   &            &                 & $0.00296$ \\
$10$&$+$&$\sqrt{2}$  &      &   &            &                 & $0.21949$ \\
 $6$&$+$&$4\sqrt{2}$ &      &   &$4\sqrt{2} $&$\simeq 5.65685$ & $3.64528$ \\
\hline
\end{tabular}
\end{center}
\end{small}
\end{table}

On the left of Fig.~\ref{fig:AB}, the averaged shelling is shown
against the radius in physical space. The result looks as erratic as
in the 1D case, except for the square-root growth of the averaged
shelling numbers at small distances. On the right, plotting the
shelling against the internal radius shows the inherent structure.
The points lie on curves that start at multiples of eight on the
ordinate and glide down convexly towards the abscissa which is met at
a point between $r^{}_1 = 1+\sqrt{2} \simeq 2.4142$ and $r^{}_2 =
(4+2\sqrt{2}\, )^{1/2} \simeq 2.6131$.  These radii correspond to the
inradius and circumradius of the window of the set $\Lambda -
\Lambda$, which is an octagon of sidelength two.

Alternatively, one may consider the corresponding random tiling. It
can be obtained, starting from the perfect pattern, by {\em simpleton
thermalization}, i.e., by repeating random simpleton flips as shown in
Fig.~\ref{fig:patch} for sufficiently many times (some 1000 flips per
vertex of the starting patch as a rule of thumb). Doing this for a
periodic approximant with 8119 vertices in the unit cell and
numerically evaluating the averaged shelling numbers for the resulting
random tiling gives the other entries of Table~\ref{tab:AB}. Note that
a real ensemble average, which would give the correct (and exact) average, is
practically impossible, so we have to rely on the self-averaging 
nature of the model.

\section{Open problems}

\hp
Illustrated by explicit results for some of the most common examples, we have
demonstrated what the averaged shelling for quasicrystals means and
how the corresponding results look like. However, closed expressions
seem difficult already for planar examples, and their real
implications are still to be understood. It is clear that several
interesting number-theoretic connections are lurking in the back, at
least if there is any reasonable generalization of the notion of theta
functions. Some further hints at this can be found in
\cite{BM}, but additional work is needed before a more definitive statement
can be made.

{}From the computational point of view, one can use the above
approach and complete the exercise for the other standard examples,
in particular for the icosahedral cases in three-dimensional space
and for the Elser-Sloane quasicrystal \cite{ES}. One
interesting aspect certainly is the possibility to distinguish perfect
from random order this way, as is apparent from the octagonal tilings
treated above. This is to be compared with the Fourier space approach
\cite{JB,JRB} and could give rise to improved statements about the
validity of perfect versus random tiling models.

\section*{Acknowledgements}

M.~B. would like to thank R.~V.~Moody and A.~Weiss for inspiring
discussions. This work was supported by the German Science Foundation
(DFG).

\vspace{3ex}
\begin{footnotesize}

\end{footnotesize}

\end{document}